\documentclass{llncs}
\usepackage{graphicx}

\usepackage{amssymb}
\usepackage{amsmath}

\newcommand{\R}{{I\!\!R}}

\def\R{{\rm I}\! {\rm R}}

\def\C{\mathbb{C}}

\usepackage[]{algorithm2e}

\begin{document}

\pagestyle{headings}

\title{Multi-stage waveform Relaxation and Multisplitting Methods for Differential Algebraic Systems}
\author{J\"urgen Geiser\inst{1} %and Fr\'{e}d\'{e}ric Magoul\`{e}s\inst{2} 
}
\institute{Ruhr University of Bochum, Department of Electrical Engineering and Information Technology, Universit\"atsstra{\ss}e 150, D-44801 Bochum, Germany, 
\\ \email{juergen.geiser@ruhr-uni-bochum.de}
% \and
%CentraleSup\'{e}lec Paris, Applied Mathematics and Systems Laboratory, High Performance Computing Research Group, Grande Voie des Vignes, 92295 Ch\^{a}tenay-Malabry Cedex, France \\
% \email{frederic.magoules@hotmail.com}
}
\maketitle

\begin{abstract}
We are motivated to solve differential algebraic equations
with new multi-stage and multisplitting methods.
The multi-stage strategy of the waveform relaxation (WR)
methods are given with outer and inner iterations.
While the outer iterations decouple the initial value problem of 
differential algebraic equations (DAEs) in the form of 
$A \frac{d y(t)}{dt} + B y(t) = f(t)$ to  $M_A  \frac{d y^{k+1}(t)}{dt} + M_1 y^{k+1}(t) = N_1 y^k(t) + N_A  \frac{d y^{k}(t)} + f(t)$, where $A = M_A - N_A$, $B = M_1 - N_1$. The inner iterations decouple further $M_1 = M_2 - N_2$ and $M_2 = M_3 - N_3$ with additional iterative processes, such that we result to invert
simpler matrices and accelerate the solver process.
The multisplitting method use additional a decomposition of the outer iterative
process with parallel algorithms, based on the partition of unity,
such that we could improve the solver method.
We discuss the different algorithms and present a first experiment based on 
a DAE system.

\end{abstract}

Multi-splitting methods, Differential-Algebraic equations.


\section{Introduction}

We are motivated to accelerate solver methods 
for differential-algebraic equations (DAEs).
Many mathematical methods are based on such combination of
differential and algebraic equations, e.g., simulations of the power systems,
constrained mechanical systems, singular perturbations, see \cite{stott_1979}
and \cite{brenan1989}.

We start with respect to assume, that the considered partial differential-algebraic equations (PDAEs) can be semi-discretized and written into an 
initial value problem of differential algebraic equations (DAEs) in the form:
\begin{eqnarray}
A \frac{d y(t)}{dt} + B y(t) = f(t) , \; t \in [t_0, T] , 
y(t_0) = y_0 ,
\end{eqnarray}
where $A \in \C^{m \times m}$ is a singular and $B \in \C^{m \times m}$ is a
non-singular complex matrix with rank $m$ and $f(t): [t_0, T] \rightarrow \C^m$
is a sufficient smooth right hand side.

For solving such problem initial value problems, waveform-relaxation (WR) methods are developed and investigated by many authors, see \cite{leim_1991}, \cite{lerara_1982} and
\cite{walle93}.

The main idea is to decompose the partitions of large systems into 
iteratively coupled smaller subsystems and solve such subsystems independently over the integration intervals (also called time-windows), see \cite{walle93}.
For the performance of the algorithms, since recent years two-stage strategy 
is introduced in WR methods, means a first splitting in blocks for pure
parallel splitting. For each processor, we apply additional a splitting
, i.e., an inner splitting instead of a direct method, see \cite{bao_2011}.
 
We propose additional a multistage splitting, means, we assume to 
split additional the inner splitting such that we could also preform
the inner splitting with parallel splitting, while only the last inner splitting
is done serial.

The new class of multi-stage waveform-relaxation (MSWR) methods are discussed 
in the following, with respect to the three-stage waveform-relaxation
(TH-SWS) method.

For simplification, we deal with
\begin{eqnarray}
\frac{d y(t)}{dt} + B y(t) = f(t) , \; t \in [t_0, T] , 
y(t_0) = y_0 ,
\end{eqnarray}
where $B = M_1 -N_1$ and the outer iteration is obtained as
\begin{eqnarray}
&& \frac{d y^{k+1}(t)}{dt} + M_1 y^{k+1}(t) = N_1 y^k(t) + f(t), \\
&& y^{k+1}(t_0) = y_0 , \; k = 1, 2, \ldots, 
\end{eqnarray}
where $B = M_1 - N_1$, then we apply the first inner iteration:
\begin{eqnarray}
&& \frac{d z^{\nu+1}(t)}{dt} + M_2 z^{\nu+1}(t) =  N_2 z^{\nu}(t)+  N_1 y^k(t) + f(t), \\
&& z^{\nu+1}(t_0) = y^k_0 , \; k = 1, 2, \ldots, \nu = 0, \ldots, \nu_k -1 ,
\end{eqnarray}
where $M_1 = M_2 - N_2$ and we obtain $y^{k+1} = z^{\nu_k}$.

The last or second inner iteration is given as:
\begin{eqnarray}
&& \frac{d \tilde{z}^{\mu+1}(t)}{dt} + M_3 \tilde{z}^{\mu+1}(t) =  N_3 \tilde{z}^{\mu}(t)+ N_2 z^{\nu}(t)+  N_1 y^k(t) + f(t), \\
&& \tilde{z}^{\mu+1}(t_0) =  z^{\nu+1}(t_0) = y^k_0 , \\
&& \; k = 1, 2, \ldots, \nu = 0, \ldots, \nu_k - 1 , \mu = 0, \ldots,  \mu_{\nu_k} - 1 , \nonumber 
\end{eqnarray}
where $M_2 = M_3 - N_3$ and we obtain $y^{k+1} = z^{\nu_k} =  \tilde{z}^{\mu_{\nu_k}}$ .
We have at least $\mu_{\nu_k}$ inner first iterations within the 
inner second iteration ${\nu_k}$ and within the $k+1$ outer iterations.
Means we deal with three-stage iterative methods.

For more flexibility in the approaches, we apply a multisplitting method,
which is based on the partition of the unit, see \cite{leary1985}.

We decompose the solution into several units, means we have:
\begin{eqnarray}
\label{synch}
&& y^{k+1} = \sum_{p = 1}^L E_{p} y^{p, k+1} , \\
&& \sum_{p = 1}^L E_{p} = I , 
\end{eqnarray}
where $I \in \C^{m \times m}$ is the unit matrix and $E_p$ are diagonal and the diagonal entries are given as $E_{p, ii} \ge 0$ for
$i = 1, \ldots, m$. 

We can solve $L$ independent waveform relaxation schemes in parallel,
while the synchronization or update is done with the Equation (\ref{synch}).

The paper is outlined as following.
In Section \ref{solver}, we discuss the different hierarchy of solver methods.
The numerical experiments are presented in Section \ref{num}.
The conclusions are done in Section \ref{conc}.

\section{Hierarchy of Solver schemes}
\label{solver}

In the following, we deal with the following hierarchy of solver schemes:
\begin{enumerate}
\item One-stage WR schemes,
\item Two-stage WR schemes,
\item Three-stage WR schemes,
\item Multisplitting WR schemes: Jacobian-, Gauss-Seidel-Types,
\end{enumerate}
where, we simplifiy the inversion of the matrices between the one-stage to the
three-stage method, means the simplification of the inversion is done via 
additional inner iteratvie stages.
Further the multisplitting approach allows to be more flexible in the
parallelization of a multi-stage method.

\subsection{One-stage WR method}

In the following, we discuss the one-stage WR method.
\begin{enumerate}
\item We have the following WR method (in parallel):
\begin{eqnarray}
\label{mobile1_1_2}
&& ( M_A + h M_1) y_{n+1}^{k+1} = h ( N_1 + \frac{1}{h} N_A) y_{n+1}^k + M_A y_n^{k+1} \\
&&  - N_A y_n^{k} + h f_{n+1} , \nonumber \\
&& y_0^{k+1} = y_0, \; k = 0, 1, \ldots, K, \; n = 0, 1,2, \ldots, J , \nonumber
\end{eqnarray}

\item
We have the following WR method (in serial):
\begin{eqnarray}
\label{mobile1_1_2}
&& ( M_A + h M_1) y_{n+1}^{k+1} = h ( N_1 + \frac{1}{h} N_A) y_{n+1}^k + M_A y_n \\
&&  - N_A y_n + h f_{n+1} , \nonumber \\
&& y_{n+1}^{0} = y_n , \; k = 0, 1, \ldots, K, \; n = 0, 1,2, \ldots, J , \nonumber 
\end{eqnarray}
where we apply the algorithm with $p = 50, q =6 , J = 20, h=0.1$, here we have $\Delta x =1.0$ and $\frac{D}{\Delta x^2} = 1.0$.

We have 2 possible stopping criteria:
\begin{enumerate}
\item Error bound: \\ 
We have an stopping error norm with $|| y_{n+1}^{k+1} - y_{n+1}^k|| \le 10^{-3}$.
\item Fix number of outer-iterative steps: $K =20 $.
\end{enumerate}

\end{enumerate}

\subsection{Two-Stage WR method}

The two-stage WR method is given as:

\begin{enumerate}
\item We have the following Two-Stage WR method (in parallel):

\begin{eqnarray}
\label{mobile1_1_2}
&& ( M_A + h M_2) z_{n+1}^{\nu+1} = h N_2 z_{n+1}^{\nu} + h ( N_1 + \frac{1}{h} N_A) y_{n+1}^k + M_A z_n^{\nu+1} \\
&& - N_A y_n^{k} + h f_{n+1} , \nonumber \\
&& z_0^{\nu+1}(t_0) = y^{k}(t_0)= y_0, \; k = 0, 1, \ldots, K, \nonumber \\
&&  \nu = 0, 1, \ldots, \nu_k, \; n = 0, 1,2, \ldots, J. \nonumber
\end{eqnarray}
\item We have the following Two-Stage WR method (in serial):

\begin{eqnarray}
\label{mobile1_1_2}
&& ( M_A + h M_2) z_{n+1}^{\nu+1} = h N_2 z_{n+1}^{\nu} + h ( N_1 + \frac{1}{h} N_A) y_{n+1}^k + M_A y_n \\
&& - N_A y_n + h f_{n+1} , \nonumber \\
&& z_{n+1}^{0}(t_n) = y_{n+1}^{k} , \; \nu = 0, 1, \ldots, \nu_k, \mbox{inner iteration}\nonumber \\
&& y_{n+1}^{k+1} = z_{n+1}^{\nu_k}, \; k = 0, 1, \ldots, K, \mbox{outer iteration} \nonumber\\
&& z_{n+1}^{0}(t_n) = y_n, \; \mbox{initialization}\nonumber \\
&&  n = 0, 1,2, \ldots, J ,\nonumber.
\end{eqnarray}

Two-Stage WR algorithm (serial) \ref{multi_2} is the given as:

\begin{algorithm}[H]
 %\KwData{this text}
 %\KwResult{how to write algorithm with \LaTeX2e }
Given the initial vector $y_0 = y(0)$ , \\
$z_{n+1}^{0}(0) = y_0,$ \\
\For
{$n = 0, 1, \ldots, J$}
{
\For
{$k = 0, 1, \ldots, K$}
{
\For
{$\nu = 0, 1, \ldots, \nu_k$}
{ $ ( M_A + h M_2) z_{n+1}^{\nu+1} = h N_2 z_{n+1}^{\nu} + h ( N_1 + \frac{1}{h} N_A) y_{n+1}^k + M_A y_n  - N_A y_n + h f_{n+1} , $ }
$ y_{n+1}^{k+1} = z_{n+1}^{\nu_k}, $ \\
$ z_{n+1}^{0} = y_{n+1}^{k+1} , $
}
$y_{n+1} = z_{n+1}^{\nu_K} , $\\
$z_{n+1}^{0}(t_{n+1}) = y_{n+1},$
}
 \caption{\label{multi_2} Two-Stage WR algorithm (serial)}
\end{algorithm}
where we apply the algorithm with $p = 50, q =6, J = 20, h=0.1$. Here we have $\Delta x =1.0$ and $\frac{D}{\Delta x^2} = 1.0$.

We have 2 possible stopping criteria:
\begin{enumerate}
\item Error bound: \\ 
We have an stopping error norm : \\
for the outer iteration with $|| y_{n+1}^{k+1} - y_{n+1}^k|| \le 10^{-3}$, \\
and for the inner iteration $|| z_{n+1}^{\nu+1} - z_{n+1}^{\nu}|| \le 10^{-3}$ 
\item Fix number of outer-iterative steps: $K = 5$ and inner iterative steps $\nu_k = 4$.
\end{enumerate}

\end{enumerate}

\subsection{Three-Stage WR method}

The three-stage WR method is given as:

\begin{enumerate}
\item We have the following Three-Stage WR method (in parallel):

\begin{eqnarray}
\label{mobile1_1_2}
&& ( M_A + h M_3) \tilde{z}_{n+1}^{\mu+1} = h N_3 \tilde{z}_{n+1}^{\mu} + h N_2 z_{n+1}^{\nu} + h ( N_1 + \frac{1}{h} N_A) y_{n+1}^k \nonumber \\
&& + M_A \tilde{z}_n^{\mu+1} - N_A y_n^{k} + h f_{n+1} , \nonumber \\
&& \tilde{z}_0^{\mu+1}(t_0) = z_0^{\nu+1}(t_0) = y^{k}(t_0)= y_0, \; k = 0, 1, \ldots, K, \nonumber \\
&&  \nu = 0, 1, \ldots, \nu_k,  \mu = 0, 1, \ldots, \mu_{\nu_k},  \; n = 0, 1,2, \ldots, J. \nonumber
\end{eqnarray}
\item We have the following Two-Stage WR method (in serial):

\begin{eqnarray}
\label{mobile1_1_2}
&& ( M_A + h M_3) \tilde{z}_{n+1}^{\mu+1} = h N_3 \tilde{z}_{n+1}^{\mu} + h N_2 z_{n+1}^{\nu} + h ( N_1 + \frac{1}{h} N_A) y_{n+1}^k \nonumber \\
&& + M_A y_n - N_A y_n + h f_{n+1} , \nonumber \\
&& \tilde{z}_{n+1}^{0}(t_n) = z_{n+1}^{\nu_k} , \; \mu = 0, 1, \ldots, \mu_{\nu_k}, \mbox{first inner iteration}\nonumber \\
&& z_{n+1}^{0}(t_n) = y_{n+1}^{k} , \; \nu = 0, 1, \ldots, \nu_k, \mbox{second inner iteration}\nonumber \\
&& y_{n+1}^{k+1} = \tilde{z}_{n+1}^{\mu_{\nu_k}}, \; k = 0, 1, \ldots, K, \mbox{outer iteration} \nonumber\\
&& \tilde{z}_{n+1}^{0}(t_n) = y_n, \; \mbox{initialization}\nonumber \\
&&  n = 0, 1,2, \ldots, J ,\nonumber.
\end{eqnarray}

Three-Stage WR algorithm (serial) \ref{multi_3} is the given as:

\begin{algorithm}[H]
 %\KwData{this text}
 %\KwResult{how to write algorithm with \LaTeX2e }
Given the initial vector $y_0 = y(0)$ , \\
$z_{n+1}^{0}(0) = y_0,$ \\
\For
{$n = 0, 1, \ldots, J$}
{
\For
{$k = 0, 1, \ldots, K$}
{
\For
{$\nu = 0, 1, \ldots, \nu_k$}
{
\For
{$\mu = 0, 1, \ldots, \mu_{\nu_k}$}
{
$( M_A + h M_3) \tilde{z}_{n+1}^{\mu+1} = h N_3 \tilde{z}_{n+1}^{\mu} + h N_2 z_{n+1}^{\nu} + h ( N_1 + \frac{1}{h} N_A) y_{n+1}^k + M_A y_n - N_A y_n + h f_{n+1} , $ }
$ z_{n+1}^{\nu+1} = \tilde{z}_{n+1}^{\mu_{\nu}}, $ \\
$ \tilde{z}_{n+1}^{0} = z_{n+1}^{\nu+1} , $
}
$ y_{n+1}^{k+1} = z_{n+1}^{\nu_k}, $ \\
$ z_{n+1}^{0} = y_{n+1}^{k+1} , $
}
$y_{n+1} = z_{n+1}^{\nu_K} , $\\
$z_{n+1}^{0}(t_{n+1}) = y_{n+1},$
}
 \caption{\label{multi_3} Three-Stage WR algorithm (serial)}
\end{algorithm}
where we apply the algorithm with $p = 50, q =6 , J = 20, h=0.1$. Further we have $\Delta x =1.0$ and $\frac{D}{\Delta x^2} = 1.0$.

We have 2 possible stopping criteria:
\begin{enumerate}
\item Error bound: \\ 
We have an stopping error norm : \\
for the outer iteration with $|| y_{n+1}^{k+1} - y_{n+1}^k|| \le 10^{-3}$, \\
and for the inner iteration $|| z_{n+1}^{\nu+1} - z_{n+1}^{\nu}|| \le 10^{-3}$ \\
and for the second inner iteration $|| \tilde{z}_{n+1}^{\mu+1} - \tilde{z}_{n+1}^{\mu}|| \le 10^{-3}$ \\ 
\item Fix number of outer-iterative steps: $K = 5$ and inner iterative steps $\nu_k = 2$, $\mu_{\nu_k} = 2$.
\end{enumerate}

\end{enumerate}

\subsection{Multisplitting WR method}

We have the following Multi-splitting WR method (in serial/parallel):
\begin{eqnarray}
\label{mobile1_1_2}
&& ( M_{A_l} + h M_{1, l}) y_{n+1}^{l, k+1} = h ( N_{1, l} + \frac{1}{h} N_{A_l}) \left( \sum_{m=1}^L E_{l, m} y_{n+1}^{l, k} \right) + M_A y_n \\
&&  - N_A y_n + h f_{n+1} , \nonumber \\
&& y_{n+1}^{l, 0} = y_n , \; l = 1, \ldots, L, \; k = 0, 1, \ldots, K, \; n = 0, 1,2, \ldots, J , \nonumber 
\end{eqnarray}
where we apply the algorithm with $p = 50, q =6 , J = 20, h=0.1$ and
the error norm $|| y_{n+1}^{k+1} - y_{n+1}^k|| \le 10^{-3}$, here we have $\Delta x =1.0$ and $\frac{D}{\Delta x^2} = 1.0$. Further $L$ are the number of the processors.

Without loosing the generality of the method, we concentrate on the
following to $L = 2$.

\subsubsection{Jacobian-Method}

The first processor is computing:
\begin{eqnarray}
\label{mobile1_1_2}
&& ( M_{A_1} + h M_{1, 1}) y_{n+1}^{1, k+1} = h ( N_{1, 1} + \frac{1}{h} N_{A_1}) \left( E_{1, 1} y_{n+1}^{1, k} + E_{1, 2} y_{n+1}^{2, k} \right) + M_{A} y_n \nonumber \\
&&  - N_{A} y_n + h f_{n+1} , \\
&& y_{n+1}^{1, 0} = y_n , \; l = 1, \ldots, L, \; k = 0, 1, \ldots, K, \; n = 0, 1,2, \ldots, J , \nonumber 
\end{eqnarray}

The second processor is computing:
\begin{eqnarray}
\label{mobile1_1_2}
&& ( M_{A_2} + h M_{1, 2}) y_{n+1}^{2, k+1} = h ( N_{1, 2} + \frac{1}{h} N_{A_2}) \left( E_{2, 1} y_{n+1}^{1, k} + E_{2, 2} y_{n+1}^{2, k} \right) + M_{A} y_n \nonumber \\
&&  - N_{A} y_n + h f_{n+1} ,  \\
&& y_{n+1}^{2, 0} = y_n , \; l = 1, \ldots, L, \; k = 0, 1, \ldots, K, \; n = 0, 1,2, \ldots, J . \nonumber 
\end{eqnarray}
where we decide if we have to switch of the mixing means:
means if we have fulfilled:
\begin{eqnarray}
\label{mobile1_1_2}
&& || \left( E_{1, 1} y_{n+1}^{1, k} + E_{1, 2} y_{n+1}^{2, k} \right) - y_{n+1}^{1, k-1} || \le || y_{n+1}^{1, k} - y_{n+1}^{1, k-1} || \\
&& || \left( E_{2, 1} y_{n+1}^{1, k} + E_{2, 2} y_{n+1}^{2, k} \right) - y_{n+1}^{2, k-1} || \le || y_{n+1}^{2, k} - y_{n+1}^{2, k-1} ||
\end{eqnarray}
we do not switch of the mixing,
but if the mixing has a larger error we have:
\begin{eqnarray}
\label{mobile1_1_2}
&&  y_{n+1}^{k} =  y_{n+1}^{1, k} , 
\end{eqnarray}

\begin{remark}

The multisplitting is switched off, if one partial solution is much more
accurate, than the other partial solution. Then we only apply the best 
approximation.

\end{remark}

\subsubsection{Gauss-Seidel-Method (decoupled version, serial with 2 processors)}

In this version, we apply the wel-known standard Gauss-Seidel method, which has the drawback of the serial treatment with the results.

The first processor is computing:
\begin{eqnarray}
\label{mobile1_1_2}
&& ( M_{A_1} + h M_{1, 1}) y_{n+1}^{1, k+1} = h ( N_{1, 1} + \frac{1}{h} N_{A_1}) \left( E_{1, 1} y_{n+1}^{1, k} + E_{1, 2} y_{n+1}^{2, k} \right) + M_{A} y_n \nonumber \\
&&  - N_{A} y_n + h f_{n+1} , \\
&& y_{n+1}^{1, 0} = y_n , \; l = 1, \ldots, L, \; k = 0, 1, \ldots, K, \; n = 0, 1,2, \ldots, J , \nonumber 
\end{eqnarray}

\begin{remark}
Here, we can apply the result of the first processor, if we
assume, that part is much more faster to the 
second processor.
\end{remark}

The second processor is computing:
\begin{eqnarray}
\label{mobile1_1_2}
&& ( M_{A_2} + h M_{1, 2}) y_{n+1}^{2, k+1} = h ( N_{1, 2} + \frac{1}{h} N_{A_2}) \left( E_{2, 1} y_{n+1}^{1, k+1} + E_{2, 2} y_{n+1}^{2, k} \right) + M_{A} y_n \nonumber \\
&&  - N_{A} y_n + h f_{n+1} ,  \\
&& y_{n+1}^{2, 0} = y_n , \; l = 1, \ldots, L, \; k = 0, 1, \ldots, K, \; n = 0, 1,2, \ldots, J . \nonumber 
\end{eqnarray}
where we decide if we have to switch of the mixing means:
means if we have fulfilled:
\begin{eqnarray}
\label{mobile1_1_2}
&& || \left( E_{1, 1} y_{n+1}^{1, k} + E_{1, 2} y_{n+1}^{2, k} \right) - y_{n+1}^{1, k-1} || \le || y_{n+1}^{1, k} - y_{n+1}^{1, k-1} || \\
&& || \left( E_{2, 1} y_{n+1}^{1, k+1} + E_{2, 2} y_{n+1}^{2, k} \right) - y_{n+1}^{2, k-1} || \le || y_{n+1}^{2, k} - y_{n+1}^{2, k-1} ||
\end{eqnarray}
we do not switch of the mixing,
but if the mixing has a larger error we have:
\begin{eqnarray}
\label{mobile1_1_2}
&&  y_{n+1}^{k} =  y_{n+1}^{1, k} , 
\end{eqnarray}

\begin{remark}
The multisplitting is switched off, if one partial solution is much more 
accurate, than the other partial solution. Otherwise, we apply the
mixture of the results based on the multisplitting method.
\end{remark}

\subsubsection{Gauss-Seidel-Method (decoupled version)}

The first processors compute:
\begin{eqnarray}
\label{mobile1_1_2}
&& \left( ( M_{A_1} + h M_{1, 1}) -  h ( N_{1, 1} + \frac{1}{h} N_{A_1}) E_{1, 1} \right)  y_{n+1}^{1, k+1} \nonumber \\
&& = h ( N_{1, 1} + \frac{1}{h} N_{A_1}) \left( E_{1, 2} y_{n+1}^{2, k} \right) + M_{A} y_n - N_{A} y_n + h f_{n+1} , \\
&& y_{n+1}^{1, 0} = y_n , \; l = 1, \ldots, L, \; k = 0, 1, \ldots, K, \; n = 0, 1,2, \ldots, J , \nonumber
\end{eqnarray}
The second processors compute:
\begin{eqnarray}
\label{mobile1_1_2}
&& \left( ( M_{A_2} + h M_{1, 2}) - h ( N_{1, 2} + \frac{1}{h} N_{A_2}) E_{2, 2} \right) y_{n+1}^{2, k+1} \nonumber \\
&& = h ( N_{1, 2} + \frac{1}{h} N_{A_2}) \left( E_{2, 1} y_{n+1}^{1, k} \right) + M_{A} y_n - N_{A} y_n + h f_{n+1} , \\
&& y_{n+1}^{2, 0} = y_n , \; l = 1, \ldots, L, \; k = 0, 1, \ldots, K, \; n = 0, 1,2, \ldots, J . \nonumber 
\end{eqnarray}

\subsubsection{Gauss-Seidel-Method (coupled version)}

Here, we apply the coupled version of the GS method, which 
means one processor is faster with the computation and the
other processor can profit from the improved computations.

The processors compute:
\begin{eqnarray}
\label{mobile1_1_2}
&& \left( ( M_{A_1} + h M_{1, 1}) -  h ( N_{1, 1} + \frac{1}{h} N_{A_1}) E_{1, 1} \right)  y_{n+1}^{1, k+1} \nonumber \\ 
&& =  h ( N_{1, 1} + \frac{1}{h} N_{A_1}) \left( E_{1, 2} \tilde{y}_{n+1}^{2, k+1} \right) + M_{A} y_n - N_{A} y_n + h f_{n+1} , \\
&& \left( ( M_{A_2} + h M_{1, 2}) - h ( N_{1, 2} + \frac{1}{h} N_{A_2}) E_{2, 2} \right) y_{n+1}^{2, k+1} \nonumber \\
&& = h ( N_{1, 2} + \frac{1}{h} N_{A_2}) \left( E_{2, 1} \tilde{y}_{n+1}^{1, k+1} \right) + M_{A} y_n - N_{A} y_n + h f_{n+1} , \\
&& y_{n+1}^{1, 0} = y_n , \; l = 1, \ldots, L, \; k = 0, 1, \ldots, K, \; n = 0, 1,2, \ldots, J , \nonumber \\
&& y_{n+1}^{2, 0} = y_n , \; l = 1, \ldots, L, \; k = 0, 1, \ldots, K, \; n = 0, 1,2, \ldots, J . \nonumber 
\end{eqnarray}
where we have two cases:
\begin{itemize}
\item If Processor $1$ is faster than Processor $2$:
\begin{eqnarray}
\label{mobile1_1_2}
&& \tilde{y}_{n+1}^{1, k+1} = y_{n+1}^{1, k+1} , \\
&& \tilde{y}_{n+1}^{2, k+1} = y_{n+1}^{2, k} ,
\end{eqnarray}
\item If Processor $2$ is faster than Processor $1$:
\begin{eqnarray}
\label{mobile1_1_2}
&& \tilde{y}_{n+1}^{1, k+1} = y_{n+1}^{1, k} , \\
&& \tilde{y}_{n+1}^{2, k+1} = y_{n+1}^{2, k+1} .
\end{eqnarray}
\end{itemize}

\begin{remark}
For all multisplitting methods, we can also extend the one-stage waveform-relaxation method to a multi-stage waveform-relaxation method.
\end{remark}

\section{Numerical Experiments}
\label{num}

In a first experiment, we apply a
partial differential algebraic equation (PDAE),
which combines partial differential and algebraic
equations.

We choose an experiment, which is based on the following
two equations:
\begin{eqnarray}
\label{mobile1_1_2}
&& \partial_t c_1  + \nabla \cdot {\bf F} c_1  =  f_1(t) ,  \; \mbox{in} \; \Omega \times [0, t] , \\
&&  \nabla \cdot {\bf F} c_2  =  f_2(t) ,  \; \mbox{in} \; \Omega \times [0, t] , \\
&& {\bf F}  = - D \nabla , 
\end{eqnarray}

and we have the following DAE problem:
\begin{eqnarray}
\label{mobile1_1_2}
&& A \partial_t c  + B c  =  f(t) ,  \; \mbox{in} \; [0, t] , 
\end{eqnarray}

The analytical solution is given as:
\begin{eqnarray}
\label{mobile1_1_2}
&& y = [\cos(t), \sin(t), t, \cos(t), \sin(t), t, \ldots, \cos(t), \sin(t), t] \in \R^m , 
\end{eqnarray}
and we have to calculate $f(t)$ as:
\begin{eqnarray}
\label{mobile1_1_2}
&& f(t) = A \partial_t y(t)  + B y(t) ,  \; \mbox{in} \; [0, t] , 
\end{eqnarray}
where $y$ is the analytical solution.

We apply the error in $L_2$ or $L_{max}$-norm means:
\begin{eqnarray}
\label{mobile1_1_2}
&& err_{L_2}(t) = \frac{1}{\Delta x} ( \sum_{i=1}^I ( y_{ana}(x_i, t) - y_{num}(x_i, t) )^2 )^{1/2} , \\ 
&& err_{max}(t) =  \max_{i=1}^I || y_{ana}(x_i, t) - y_{num}(x_i, t) || . 
\end{eqnarray}

In the following we deal with the semidiscretized equation given with the 
matrices:
\begin{eqnarray}
\label{eq20}
&& A =  
\left(
\begin{array}{c c c c c}
 I &  & & &  \\
  & \ddots & & &     \\
  &  & I & &     \\
   &     &   & 0 &    \\
   &    &  &  & 0
\end{array}
\right)  \in \R^{m \times m} , 
\end{eqnarray}
where $I, 0 \in \R^{p \times p}$

We have the following two operators for the splitting method:
\begin{eqnarray}
B_1 & = & \left(\begin{array}{rrrrr}
 4 & -1 & ~ & ~ & ~ \\
 -1 & 4 & -1 & ~ & ~ \\
 ~ & \ddots & \ddots & \ddots & ~ \\
 ~ & ~ & -1 & 4 & -1 \\
 ~ & ~ & ~ & -1 & 4
\end{array}\right) \in \R^{p \times p}
\end{eqnarray}
\begin{eqnarray}
B & = &  \frac{D}{\Delta x^2}\cdot  \left(\begin{array}{rrrrr}
 B_1 & -I & ~ & ~ & ~ \\
 -I & B_1 & -I & ~ & ~ \\
 ~ & \ddots & \ddots & \ddots & ~ \\
 ~ & ~ & -I & B_1 & -I \\
 ~ & ~ & ~ & -I & B_1
\end{array}\right) \in \R^{m\times m}
\end{eqnarray}
with $p q = m$, where we assume $ \frac{D}{\Delta x^2} = 1$.

Means $A, B$ are $m \times m$ block-matrices.

We have the following splitting:
\begin{eqnarray}
\label{eq20}
&& N_A = 
\left(
\begin{array}{c c c c c}
 0 &  & & &  \\
  & \ddots & & &     \\
  &  & 0 & &     \\
   &     &   & \frac{1}{100} I  &    \\
   &    &  &  & 0
\end{array}
\right) , 
\end{eqnarray}
\begin{eqnarray}
N_1 & = &  \frac{D}{\Delta x^2}\cdot \left(\begin{array}{rrrrr}
 2 I & I & ~ & ~ & ~ \\
 I & 2 I & I & ~ & ~ \\
 ~ & \ddots & \ddots & \ddots & ~ \\
 ~ & ~ & I & 2 I & I \\
 ~ & ~ & ~ & I & 2 I
\end{array}\right) \in \R^{m \times m}
\end{eqnarray}
\begin{eqnarray}
M_2 & = & \frac{D}{\Delta x^2}\cdot \left(\begin{array}{rrrrr}
 8 I & 0 & ~ & ~ & ~ \\
 0 & 8 I & I & ~ & ~ \\
 ~ & \ddots & \ddots & \ddots & ~ \\
 ~ & ~ & 0 & 8 I & 0 \\
 ~ & ~ & ~ & 0 & 8 I
\end{array}\right) \in \R^{m \times m}
\end{eqnarray}
where $I, 0 \in \R^{p \times p}$
\begin{eqnarray}
M_3 & = & \frac{D}{\Delta x^2}\cdot \left(\begin{array}{rrrrr}
 10 I & 0 & ~ & ~ & ~ \\
 0 & 10 I & I & ~ & ~ \\
 ~ & \ddots & \ddots & \ddots & ~ \\
 ~ & ~ & 0 & 10 I & 0 \\
 ~ & ~ & ~ & 0 & 10 I
\end{array}\right) \in \R^{m \times m}
\end{eqnarray}
where $I, 0 \in \R^{p \times p}$

We have the following operators:
\begin{eqnarray}
\label{mobile1_1_2}
&& M_A = A +  N_ A , \\
&& M_1 = B + N_1  , \\
&& N_2 = M_2 - M_1 , \\
&& N_3 = M_3 - M_2 , 
\end{eqnarray}

Further we have the following matrices:
\begin{eqnarray}
\label{eq20}
&& N_{A_1} = 
\left(
\begin{array}{c c c c c}
 0 &  & & &  \\
  & \ddots & & &     \\
  &  & 0 & &     \\
   &     &   & \frac{1}{100} I  &    \\
   &    &  &  & 0
\end{array}
\right) , 
 N_{A_2} = 
\left(
\begin{array}{c c c c c}
 0 &  & & &  \\
  & \ddots & & &     \\
  &  & 0 & &     \\
   &     &   & 0   &    \\
   &    &  &  & \frac{1}{100} I
\end{array}
\right) \in \R^{m \times m} , 
\end{eqnarray}
\begin{eqnarray}
&& N_{1,1} = \frac{D}{\Delta x^2}\cdot \left(\begin{array}{rrrrr}
 2 I & I & ~ & ~ & ~ \\
 I & 2 I & I & ~ & ~ \\
 ~ & \ddots & \ddots & \ddots & ~ \\
 ~ & ~ & I & 2 I & I \\
 ~ & ~ & ~ & I & 2 I
\end{array}\right)  \in \R^{m \times m} , \\
&& N_{1,2} = \frac{D}{\Delta x^2}\cdot \left(\begin{array}{rrrrr}
 3 I & I & ~ & ~ & ~ \\
 I & 3 I & I & ~ & ~ \\
 ~ & \ddots & \ddots & \ddots & ~ \\
 ~ & ~ & I & 3 I & I \\
 ~ & ~ & ~ & I & 3 I
\end{array}\right) \in \R^{m \times m}
\end{eqnarray}
and the overlapping matrices are given as (where we deal with
a symmetric overlap, means in both directions of the diagonal matrix):
\begin{enumerate} 
\item Overlap is between block-matrices at $m/2$ and $m/2 + 1$):
\begin{eqnarray}
\label{eq20}
&& E_{1, 1} = 
\left(
\begin{array}{c c c c c c c}
 I &  & & & & & \\
  & \ddots & & & & &    \\
  &  & I & &  & &  \\
   &     &   & \alpha_1 I  & & &  \\
   &     &   &   & 0 & &  \\
   &    &  &  &  & \ddots& \\
   &    &  &  &  &  & 0
\end{array}
\right) ,  \\
&& \nonumber \\
&&
 E_{1, 2} = 
\left(
\begin{array}{c c c c c c c}
 0 &  & & & & & \\
  & \ddots & & & & &    \\
  &  & 0 & &  & &  \\
   &     &   & \alpha_2 I  & & &  \\
   &     &   &   & I & &  \\
   &    &  &  &  & \ddots& \\
   &    &  &  &  &  & I
\end{array}
\right) \in \R^{m \times m} , 
\end{eqnarray}
and the next decomposition:
\begin{eqnarray}
\label{eq20}
&& E_{2, 1} = 
\left(
\begin{array}{c c c c c c c}
 I &  & & & & & \\
  & \ddots & & & & &    \\
  &  & I & &  & &  \\
   &     &   & \alpha_3 I  & & &  \\
   &     &   &   & 0 & &  \\
   &    &  &  &  & \ddots& \\
   &    &  &  &  &  & 0
\end{array}
\right) ,  \\
&& \nonumber \\
&&
 E_{2, 2} = 
\left(
\begin{array}{c c c c c c c}
 0 &  & & & & & \\
  & \ddots & & & & &    \\
  &  & 0 & &  & &  \\
   &     &   & \alpha_4 I  & & &  \\
   &     &   &   & I & &  \\
   &    &  &  &  & \ddots& \\
   &    &  &  &  &  & I
\end{array}
\right) \in \R^{m \times m} , 
\end{eqnarray}
where $\alpha_1 + \alpha_2 = 1$ and $\alpha_3 + \alpha_4 = 1$,
Here we have the overlap $o= 1$, means $E_{1,1}$ and $E_{2,2}$ have only one
line overlap with $I$. 
\begin{remark}
An extension is to apply different overlapping areas in the $E_{1}$ and
$E_2$ decomposition. 
\end{remark}

\item The largest overlap is $o = m/2 -1$ (where we assume $m$ is even), means we overlap nearly the full matrices except the lowest and uppermost entry,
see:
\begin{eqnarray}
\label{eq20}
&& E_{1, 1} = 
\left(
\begin{array}{c c c c c c c c}
 I &    &   & & & & & \\
   &  \alpha_1 I &   & & & & &    \\
   &    & \ddots &   &   &   &   &     \\
   &    &        & \alpha_1 I  &  &   &  &  \\
   &    &        &    & \alpha_1 I  &         &  &  \\
   &    &        &    &    & \ddots  &  &   \\
   &    &        &    &    &  & \alpha_1 I &  \\
   &    &        &    &     &     &    & 0
\end{array}
\right) ,  \\
&& \nonumber \\
&&
 E_{1, 2} = 
\left(
\begin{array}{c c c c c c c c}
 0 &    &   & & & & & \\
   &  \alpha_2 I &   & & & & &    \\
   &    & \ddots &   &   &   &   &     \\
   &    &        & \alpha_2 I  &  &   &  &  \\
   &    &        &    & \alpha_2 I  &         &  &  \\
   &    &        &    &    & \ddots  &  &   \\
   &    &        &    &    &  & \alpha_2 I &  \\
   &    &        &    &     &     &    & I
\end{array}
\right) \in \R^{m \times m} , 
\end{eqnarray}
and the next decomposition:
\begin{eqnarray}
\label{eq20}
&& E_{2, 1} = 
\left(
\begin{array}{c c c c c c c c}
 I &    &   & & & & & \\
   &  \alpha_3 I &   & & & & &    \\
   &    & \ddots &   &   &   &   &     \\
   &    &        & \alpha_3 I  &  &   &  &  \\
   &    &        &    & \alpha_3 I  &         &  &  \\
   &    &        &    &    & \ddots  &  &   \\
   &    &        &    &    &  & \alpha_3 I &  \\
   &    &        &    &     &     &    & 0
\end{array}
\right) , \\
&& \nonumber \\
&&
 E_{2, 2} = 
\left(
\begin{array}{c c c c c c c c}
 0 &    &   & & & & & \\
   &  \alpha_4 I &   & & & & &    \\
   &    & \ddots &   &   &   &   &     \\
   &    &        & \alpha_4 I  &  &   &  &  \\
   &    &        &    & \alpha_4 I  &         &  &  \\
   &    &        &    &    & \ddots  &  &   \\
   &    &        &    &    &  & \alpha_4 I &  \\
   &    &        &    &     &     &    & I
\end{array}
\right) \in \R^{m \times m} .
\end{eqnarray}
\begin{remark}
An extension is to apply different overlapping areas in the $E_{1}$ and
$E_2$ decomposition. 
\end{remark}

\end{enumerate}

Further we have
\begin{eqnarray}
\label{eq20}
 E_1 = E_2 = 
\left(
\begin{array}{c c c c c c c}
 I &  & & & & & \\
  & \ddots & & & & &    \\
  &  & I & &  & &  \\
   &     &   & I  & & &  \\
   &     &   &   & I & &  \\
   &    &  &  &  & \ddots& \\
   &    &  &  &  &  & I
\end{array}
\right) \in \R^{m \times m}, 
\end{eqnarray}
We have the following operators:
\begin{eqnarray}
\label{mobile1_1_2}
&& M_{A_1} = A +  N_{A_1} , \\
&& M_{1,1} = B + N_{1,1}  , \\
&&  M_{A_2} = A +  N_{A_2} , \\
&& M_{1,2} = B + N_{1,2}  , \\
&& E_1 = E_{1,1} + E_{1,2} , \\
&& E_2 = E_{2,1} + E_{2,2} ,
\end{eqnarray}

\begin{remark}
We compared the errors between the multi-level WR and the Multisplitting WR with Jacobian and Gauss-Seidel
types.
For the multi-level methods, we present the benefit in the higher level methods, while we only
invert smaller matrices. The highest accuracy is given with the one-level method and the
MS-Gauss-Seidel method, while the inversion matrix has the largest amount of information, but the
methods are at least very expensive.

We could also improve the accuracy of the MS methods based on the different
overlapping  means for $o=1$, we have only one overlap, while $o=m/2-1$ has the largest overlap.
The balance and optimal values are between.

\end{remark}

\begin{figure}[ht]
\begin{center}  
\includegraphics[width=5.0cm,angle=-0]{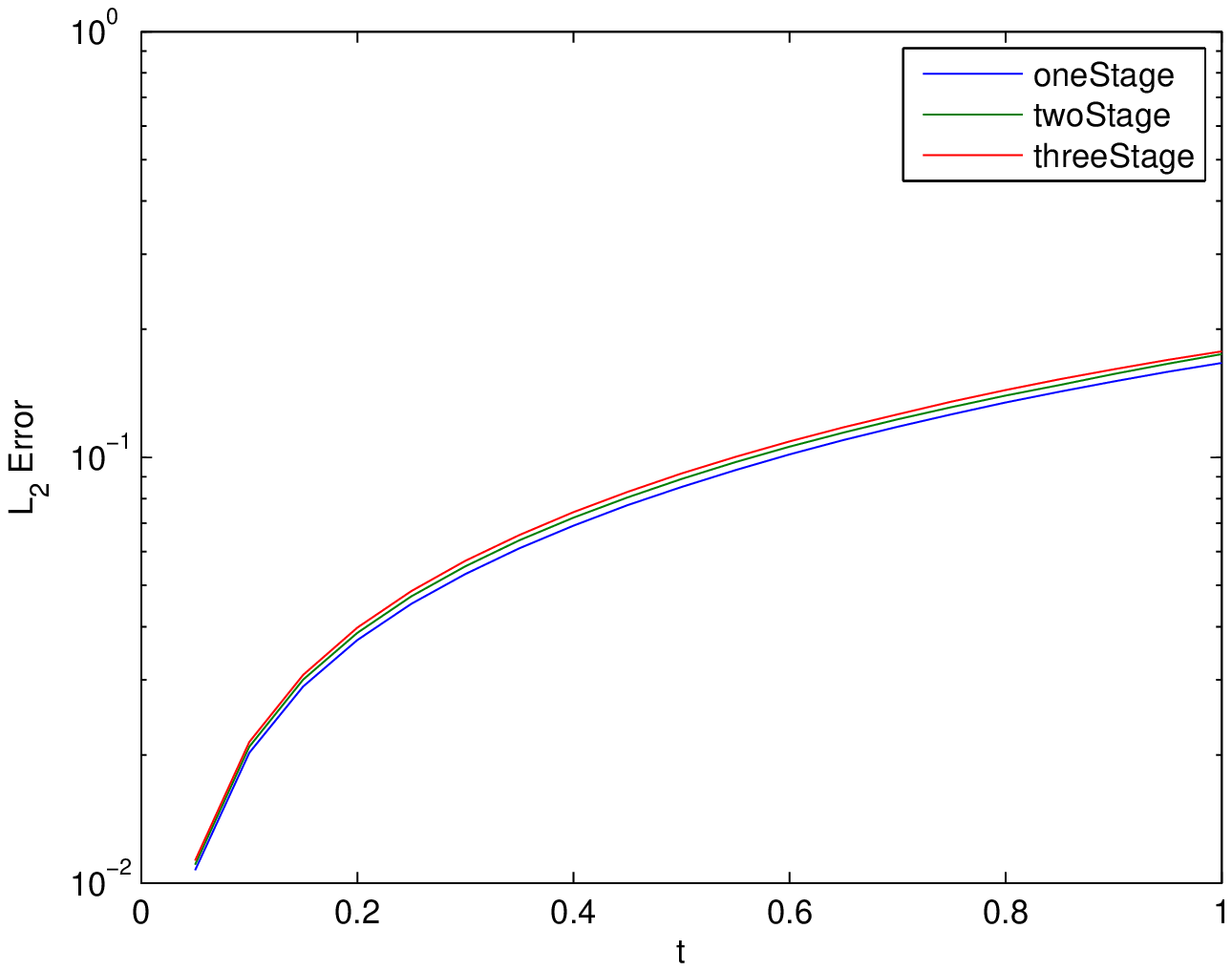}
\includegraphics[width=5.0cm,angle=-0]{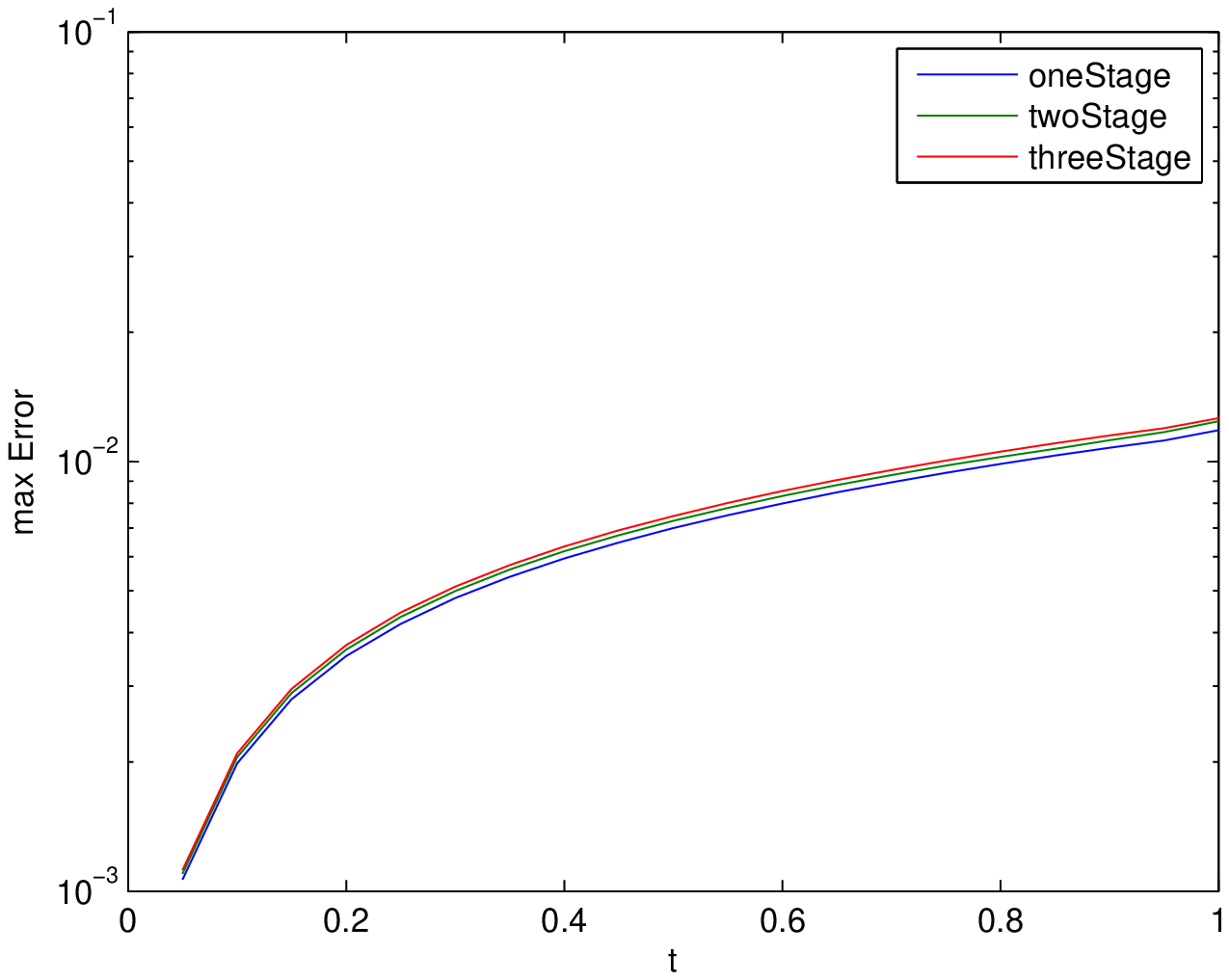}
\end{center}
\caption{\label{one-three-stage} The errors between the exact and numerical scheme of the one-, two- and three-stage method is given (left hand side: $L_2$-errors, right hand side: $L_{\infty}$-error).}
\end{figure}

We apply the numerical example and obtain
for the one-stage, two-stage and three-stage method the
following results in Figure \ref{one-three-stage}.

\begin{remark}
We obtain the same accuracy of the three and two-stage method as
for the one-stage method. This means, that we can reduce the computational 
amount of work and received the same accurate result.
\end{remark}

\begin{figure}[ht]
\begin{center}  
\includegraphics[width=5.0cm,angle=-0]{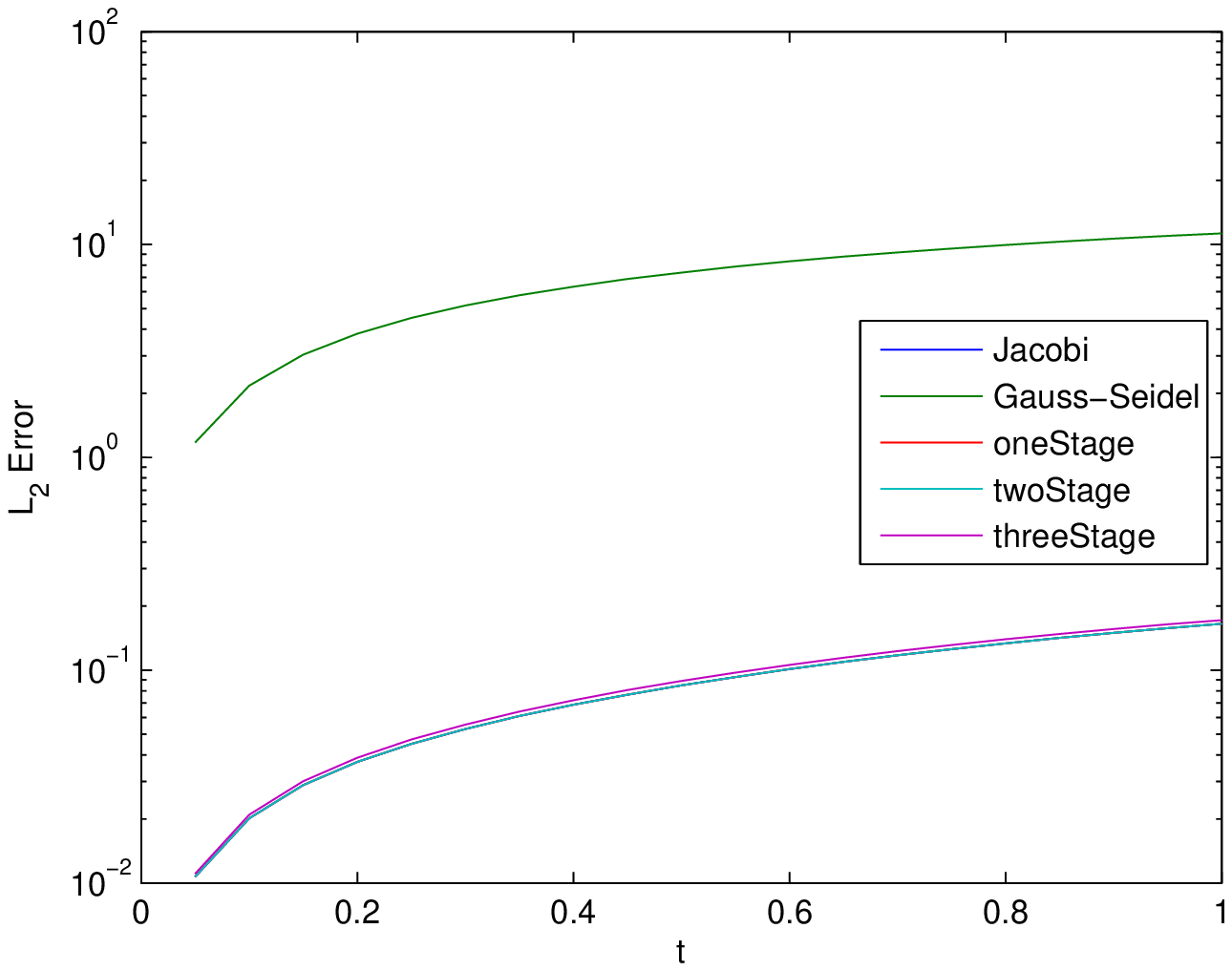}
\includegraphics[width=5.0cm,angle=-0]{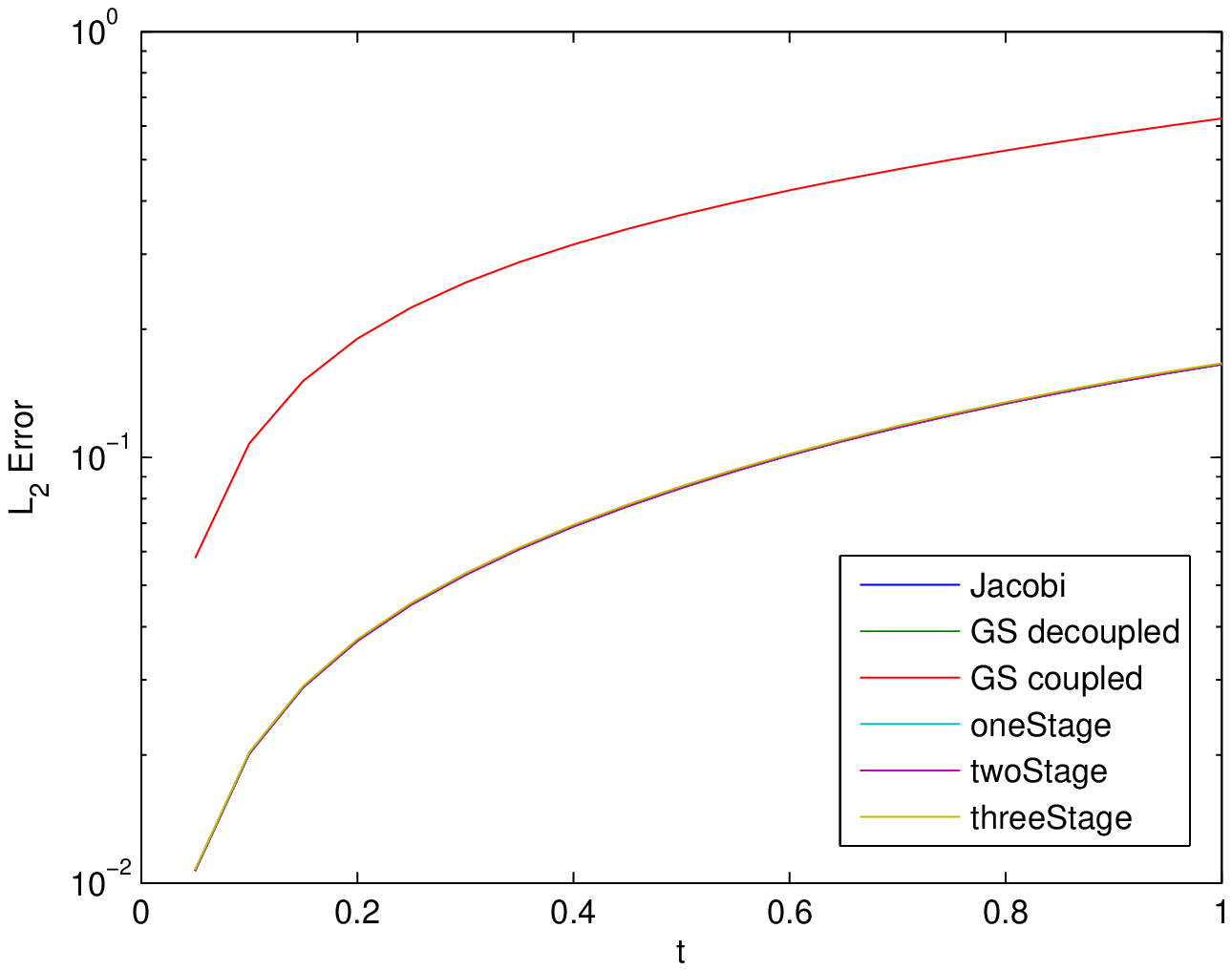}
\end{center}
\caption{\label{multi_1} The errors between the exact and numerical scheme of the one-stage, two-stage, three-stage, Jacobi- and GS Multisplitting method is presented (left hand side: one-stage, two-stage, three-stage, Jacobi-, uncoupled GS-method, right hand side: one-stage, two-stage, three-stage, Jacobi-, uncoupled and coupled GS-methods.}
\end{figure}

We apply the numerical example and obtain
for the Multisplitting method the
following results in Figure \ref{multi_1}.

\begin{remark}
We obtain the same accuracy of the three and two-stage method as
for the one-stage method. This means, that we can reduce the computational 
amount of work and received the same accurate result. The multisplitting
method has also the same accuracy as the different multi-stage methods, here,
we have the benefit of the parallel versions.
\end{remark}

\section{Conclusions and Discussions}
\label{conc}

We discuss multi-stage waveform-relaxation methods and
multisplitting methods for differential algebraic equations.
While the multi-stage waveform-relaxation methods can reduce their
computational work with simplifying the inverse matrices, the
multisplitting methods have their benefits in parallelizing their
procedure. We test the ideas in a first partial differential
algebraic equation and see the benefit in the multi-stage waveform-relaxation
method.
In future, we will discuss the numerical analysis of the different
methods and present more numerical examples.

\bibliographystyle{plain}

\end{document}